
\magnification=1200

\font\title=cmr10 scaled\magstep2
\font\small=cmr7 scaled\magstep1
\font\smallb=cmb8 scaled\magstep1

\noindent
{\title Reduction of Linear Programming to Linear Approximation}
\bigskip

\noindent
LEONID N. VASERSTEIN, 

\noindent
{\it Department of Mathematics, Penn State U., University Park, PA 16802 }

\noindent
{\it 
(e-mail: vstein@math.psu.edu)}
\bigskip

\noindent
{\small Received  Jan. 27, 2006}
\bigskip

\noindent
{\smallb Abstract.}
{\small 
It is well known that every  $l^{\infty}$ linear approximation problem can be reduced to a linear program. In this paper we show that conversely 
every linear program can be reduced to an  $l^{\infty}$ linear approximation problem.}

\bigskip

\noindent
{\smallb Key words} 
{\small Linear programming, linear approximation, Chebyshev approximation.}

\bigskip
\noindent
It is well known that every  $l^{\infty}$ linear approximation problem can be reduced to a linear program. In this paper we show that conversely 
every linear program can be reduced to an  $l^{\infty}$ linear approximation problem.

Now we recall relevent definitions. 

An {\it affine function}  of variables $x_1,\ldots,x_n$ is
$b_0+c_1x_1+\cdots + c_nx_n$ where $b_0,c_i$ are given numbers.

A {\it linear constraint}  is any of the following constraints: $f \le g, 
f \ge g,  f = g, $ where $f,g$ are affine functions.

A {\it linear program} is an optimization (maximization or minimization) of an affine function
subject to a finite system of linear constraints.

An  $l^{\infty}$  linear approximation  problem, also known as (discrete)  {\it Chebyshev 
approximation problem } or finding the least-absolute-deviation fit,
 is the problem of minimization
of  the following function:
$$\max(|f_1|,\ldots , |f_m|) = \|(f_1,\ldots,f_m)\|_{\infty},$$
where $f_i$ are affine functions.  This objective function is piece-wise linear and convex. 

Given any Chebyshev approximation problem,  here is a well-known reduction (Vaserstein, 2003)  to a linear program
with one additional variable $t$:
$$t \to \min, \  {\rm subject\  to} \  -t \le f_i  \le t \  {\rm for } \  i =1,\ldots, m.$$
This is a linear program with $n+1$ variables and $2m$ linear constraints.

Now we want to reduce an arbitrary linear program to a Chebyshev approximation
problem. First of all, it is well known (Vaserstein, 2003) that every linear program can be reduced to
solving a symmetric matrix  game.

So we start with a matrix game, with the payoff matrix $M = - M^T$ of size $N$ by $N.$
Our problem is to find   a column $x=(x_i)$ (an optimal strategy)
such that  
$$  Mx \le 0, x \ge 0, \sum x_i = 1. \eqno{(1)}$$

As usual, $x \ge 0$ means that every entry of the column $x$ is  $\ge 0.$
Later we write $y \le t$ for a column $ y $ and a number $t$ if every entry of $y$ is $\le t.$  We go even further in abusing notation, denoting
by $y-t$ the column obtaining from $y$ by subtracting $t$ from every entry.
Similarly we denote by $M+c$ the matrix obtained from $M$ by adding a number $c$ to every entry.

This problem (1) (of finding an optimal strategy)  is about finding a feasible solution for a system of linear constraints. It can be written as the following linear program with an additional variable $t$:
$$  t \to \min, Mx \le t,  x \ge 0, \sum x_i = 1.\eqno{(2)}$$

Now we find the largest entry $c$ in the matrix $M$.  If $c=0,$
then $M=0$ and the problem (1) is trivial (every mixed strategy $x$ is optimal). So we assume that $c > 0.$

Adding the number $c$  to every entry of the matrix
$M,$ we obtain a matrix $M+c \ge 0 $ (all entries $\ge 0).$  The linear program (2)
is equivalent to
$$  t \to \min, (M+c)x \le t,  x \ge 0, \sum x_i = 1 \eqno{(3)}$$
in the sense that these two programs have the same feasible solutions and  the same optimal solutions. The optimal value for (2) is 0 while the optimal value for (3) is $c.$

Now we can rewrite (3) as follows:
$$  \|(M+c)x \|_{\infty}  \to \min,  x \ge 0, \sum x_i = 1 \eqno{(4)}$$
which is a Chebyshev approximation problem with additional linear constraints.
We used that $M+c \ge 0,$ hence  $(M+c)x  \ge 0$ for every feasible solution $x$ in (2).
The optimal value is still $c.$

Now we rid off the constraints in(4)  as follows:
$$ \|\pmatrix{(M+c)x \cr c-x \cr  \sum x_i +c - 1 \cr -\sum x_i -c + 1}  \|_{\infty} \to \min. \eqno{(5)}$$

Note that the optimization problems (4) and (5) have the same optimal value
$c$ and every optimal solution  of (4) is optimal for (5).
Conversely, 
 for every $x$ with a negative entry, the objective function in (5) is
$> c.$ Also, for  every $x$ with  $\sum x_i \ne 1,$  the objective function in (5) is
$> c.$ So every optimal solution  for (5) is feasible and hence optimal for (4).

Thus, we have reduced solving any symmetric matrix game with $N \times N$ payoff matrix to a Chebyshev approximation problem (5) with $2N+2$ affine functions in $N$ variables.

\smallskip
 \noindent
{\bf  Remark. }   It is well known that every  $l^{1}$ linear approximation problem can be reduced to a linear program.
 Our result implies that every  $l^{1}$ linear approximation problem can be reduced to a
 $l^{\infty}$ linear approximation problem. I do not know whether the converse is true.

 Note that our reduction of the  $l^{1}$ linear approximation problem
 $$
 \sum _{i=1}^m |f_i| \to \min \eqno{(6)}
 $$
 where $f_i$ are affine functions in $n$ variables, produces first  the well-known linear program (Vaserstein, 2003)
 $$
 \sum _{i=1}^m t_i \to \min, -t_i \le f_i \le t_i
 $$
 with  $m+n$ variables and $2m$ linear constraints, then a symmetric
 game with  the payoff matrix of size  
 $(3m+2n+1) \times  (3m+2n+1),$ and finally a  
 Chebyshev approximation problem with   $6m+4n+4$   affine functions in   $3m+2n+1$ variables.
 
 By comparison,
 an obvious    direct reduction  produces 
 $$\max| f_1\pm f_2 \pm \cdots \pm f_m| \to \min
 $$
 which is a Chebyshev approximation problem  with  $2^{m-1}$ affine functions in $n$ variables.
 So this reduction increases the size exponentially, while our reduction increases size linearly.
 
 \bigskip
 \noindent
 {\bf References}
 \smallskip
 
 \noindent
Vaserstein, L. N. (2003),  {\it  Introduction to Linear Programming,}  Prentice Hall. 
(There is a Chinese translation by    Mechanical Industry Publishing House ISBN: 7111173295. )
\end